\documentclass[12pt ]{amsart}

\usepackage{amsmath}
\usepackage{amsfonts}
\usepackage{amssymb}
\usepackage{url}

\numberwithin{equation}{section}

\newtheorem{theorem}[subsection]{Theorem}
\newtheorem{lemma}[subsection]{Lemma}

\newcommand{\half}{\tfrac{1}{2}}
\newcommand{\f}{f\times \chi}
\newcommand{\summ}{\mathop{{\sum}^{\star}}}
\newcommand{\chiq}{\chi \bmod q}

\newcommand{\chid}{\chi \bmod d}

\newcommand{\V}{V\Big(\frac{nm}{q^2}\Big)}

\begin{document}

\baselineskip=17pt

\title{The second moment of Dirichlet twists of Hecke $L$-functions }

\author{Peng Gao}
\address{Nanyang Technological University, Division of Mathematical Sciences, 21 Nanyang Link Singapore 637371}
\email{penggao@ntu.edu.sg}

\author{Rizwanur Khan}
\address{University of California, Los Angeles,
         Department of Mathematics, 520 Portola Plaza, Los Angeles, CA 90095-1555, USA}
\email{rrkhan@math.ucla.edu}

\author{Guillaume Ricotta}
\address{Universit\'{e} de Bordeaux 1, Institut de Math\'{e}matiques de Bordeaux, Laboratoire A2X, B\^{a}timent A33, Bureau 301 M, 351 cours de la Lib\'{e}ration, 33405 Talence Cedex, France}
\email{guillaume.ricotta@math.u-bordeaux1.fr}

\maketitle

\section{Introduction}

Let $f$ be a fixed holomorphic Hecke eigenform of level 1 and
weight $k$. For $\Im(z)>0$ we have a Fourier expansion of $f$:
\begin{align*}
f(z)=\sum_{n\ge 1} a_f(n)n^{\frac{k-1}{2}}e^{2\pi i n z},
\end{align*}
where the coefficients $a_f(n)$ are real and satisfy the
Ramanujan-Petersson bound $|a(n)|\le d(n)$, and we have normalized
so that $a_f(1)=1$. Let $\chi$ be a Dirichlet character of modulus
$q$, and form the Dirichlet series of the twist of $f$ by $\chi$
given by
\begin{align*}
L(\f,s)= \sum_{n\ge 1} \frac{a_f(n)\chi(n)}{n^s}=\prod_{p} \Big(1-\frac{a_f(p)\chi(p)}{p^{s}}+\frac{\chi(p^2)}{p^{2s}}\Big)^{-1}
\end{align*}
for $\Re(s)>1$. When $\chi$ is primitive this is an $L$-function,
which continues to an entire function and satisfies the functional
equation
\begin{align}
(\tfrac{q}{2\pi})^{s}\Gamma(\tfrac{k}{2}+s)L(\f,\half+s)=\imath_{\chi}(\tfrac{q}{2\pi})^{-s}\Gamma(\tfrac{k}{2}-s)L(f\times
\overline{\chi},\half-s).
\end{align}
where $\imath_{\chi}=i^k \frac{\tau(\chi)^2}{q}$ and $\tau(\chi)$ is
the Gauss sum associated to $\chi$ (thus $|\imath_{\chi}|=1$). cf.
\cite[Prop. 14.20]{iwakow}.

For a large modulus $q$, we are interested in the second
power moment of $L(\f, \half)$ as $\chi$ ranges over all primitive
characters modulo $q$. Stefanicki \cite{stef} proved
\begin{align}
\label{stef}\summ_{\chiq} |L(\f, \half)|^2 =K P_q(1) \psi(q) q \log
q + O(2^{\nu(q)}q(\log q)^{0.935}),
\end{align}
where $\star$ denotes summation only over primitive characters, $K$
is a constant, $P_q(1) \psi(q)$ is a finite Euler product over the
prime divisors of $q$, and $\nu(q)$ is the number of distinct prime
factors of $q$. A modulus $q$ may have $\frac{\log q}{\log\log q}$
distinct prime factors, but the main term in (\ref{stef}) fails to
dominate the error term when $\nu(q)$ is larger than $\frac{1}{10}
\log\log q$. The condition $\nu(q)<\frac{1}{10} \log\log q$ holds
for very few integers $q$, a set of zero density in fact, by a
theorem of Hardy and Ramanujan \cite{hardyrama}. We are interested
in the following problem: can we establish an asymptotic for the
second moment for a larger class of integers $q$; if not all
integers, then almost all? This can be considered an analogue of a
problem studied by Soundararajan: in \cite{sound}, he adapts
Heath-Brown's \cite{hb} asymptotic formula for the fourth power
moment of Dirichlet $L$-functions to hold for all moduli.
Unfortunately his methods do not seem to immediately carry over to
this case and we must do something different. We should also mention
that it is an open problem to find an asymptotic expression with a
power saving error term for the left hand side of (\ref{stef}), even
for prime $q$. The corresponding problem for the fourth power moment
of Dirichlet $L$-functions with prime moduli was recently solved by
Young \cite{young}.

Let $K=\tfrac{6(4\pi)^k}{\pi\Gamma(k)}\|f\|^2$, where the norm is the Petersson norm. Define the multiplicative function $\psi(q)$ by setting $\psi(p)=(1-\tfrac{2}{p})$ and $\psi(p^k)=(1-\frac{1}{p})^2$ for $k\ge 2$, so that $q\psi(q)$ is the number of primitive characters modulo $q$. Let
\begin{align*}
P_q(s)=\prod_{p|q}\big(1-\tfrac{1}{p^s}\big)^2
\big(1-\tfrac{a_f(p)^2-2}{p^s}+\tfrac{1}{p^{2s}}\big)
\big(1-\tfrac{1}{p^{2s}}\big)^{-1}.
\end{align*}
Our main theorem is
\begin{theorem}\label{mainthm}
For integers $q$ satisfying
\begin{align}
\label{assumption} \sum_{\substack{p|q\\ p> x}} \frac{1}{p} \le
(\log \log q)^{-10},
\end{align}
where $x=\exp\big(\frac{\log\log q}{200 \log \log \log q}\big)$, we
have
\begin{align}
\summ_{\chiq} |L(\f, \half)|^2 = K P_q(1) \psi(q) q \log
q\Big(1+O\big((\log \log q)^{-1}\big) \Big),
\end{align}
where the implied constant depends on $f$.
\end{theorem}
\noindent Condition (\ref{assumption}) holds when for example
$\nu(q)\ll \exp\big(\frac{\log\log q}{300 \log \log \log q}\big)$,
which is certainly true of almost all integers. We will need the
following consequence of this condition.

\begin{lemma}\label{condition} If (\ref{assumption}) holds then we have
\begin{align}
\sum_{\substack{d|q\\ d\ge (\log q)^{0.05}}} \frac{|\mu(d)|
\prod_{p|d}(1+\tfrac{10}{\sqrt{p}}) }{d} \ll (\log \log q)^{-9}.
\end{align}
\end{lemma}
\proof Using `Rankin's trick' we have
\begin{align}
\sum_{\substack{d|q\\ d\ge (\log q)^{0.05}\\ p|d\Rightarrow p\le x}}
&\frac{|\mu(d)| \prod_{p|d}(1+\tfrac{10}{\sqrt{p}}) }{d}\\
&\nonumber \le \sum_{\substack{d|q\\  p|d\Rightarrow p\le
x}}\Big(\frac{d}{(\log q)^{0.05}}\Big)^{(\log x)^{-1}}
\frac{|\mu(d)| \prod_{p|d}(1+\tfrac{10}{\sqrt{p}}) }{d}\\
&\nonumber \ll (\log \log q)^{-10} \prod_{p\le x}
\Big(1+\frac{1}{p^{1-(\log x)^{-1}}}\Big)\ll (\log \log q)^{-9}.
\end{align}
Thus
\begin{align}
\label{trick}\sum_{\substack{d|q\\ d\ge (\log q)^{0.05}}}
&\frac{|\mu(d)| \prod_{p|d}(1+\tfrac{10}{\sqrt{p}}) }{d}\\
&\nonumber \ll \sum_{\substack{d|q\\ d\ge (\log q)^{0.05}\\
p|d\Rightarrow p> x}} \frac{|\mu(d)|
\prod_{p|d}(1+\tfrac{10}{\sqrt{p}}) }{d} + (\log\log
q)^{-9}\\
&\nonumber \ll \Big(\sum_{\substack{p|q\\p>x}}\frac{1}{p}\Big)
\Big(\sum_{d|q}\frac{|\mu(d)|
\prod_{p|d}(1+\tfrac{10}{\sqrt{p}})}{d}\Big)+ (\log\log q)^{-9}.
\end{align}
This completes the proof as the sum over divisors in the last line
of (\ref{trick}) is $\ll \prod_{p|q}(1+1/p)\ll \log \log q.$
\endproof

Throughout the paper all implied constants may depend implicitly on
$f$.

\section{Proof of Theorem \ref{mainthm}}

The following standard results can be found in sections 2 and 3.1
respectively of \cite{stef}.
\begin{lemma}[\bf Orthogonality]\label{ortho} For $(nm,q)=1$ we have
\begin{align}
\summ_{\chiq} \chi(n)\overline{\chi}(m) = \sum_{\substack{d|q\\
n\equiv m \bmod q/d}} \mu(d) \phi(q/d).
\end{align}
\end{lemma}
\begin{lemma}[\bf Approximate functional equation]\label{apprfunceqn}
If $\chi$ is a primitive character modulo $q$ we have
\begin{align}
|L(\f,\half)|^2 =\sum_{n,m}
\frac{a_f(n)a_f(m)\chi(n)\overline{\chi}(m)}{\sqrt{nm}}\V,
\end{align}
where $V(x)=\frac{1}{\pi i} \int_{(c)} G(y) x^{-y} \frac{dy}{y}$ for
any $c>0$ and
$G(s)=\frac{\Gamma(k/2+s)^2}{(2\pi)^{2s}\Gamma(k/2)^2}$ decays
rapidly in vertical lines. We have that $V(x)\ll_C \min \{1,
x^{-C}\}$ for any $C>0$, so that the sum above is essentially
supported on $nm< q^{2+\epsilon}$ for any $\epsilon>0$.
\end{lemma}

\noindent Thus we have
\begin{align}
\summ_{\chiq} |L(\f, \half)|^2&= \sum_{n,m} \summ_{\chiq}
\frac{a_f(n)a_f(m)}{\sqrt{nm}}\V\summ_{\chiq}\chi(n)\overline{\chi}(m)\\
\nonumber &= \sum_{d|q} \mu(d) \phi(q/d) \sum_{\substack{n\equiv m
\bmod q/d\\(mn,q)=1}} \frac{a_f(n)a_f(m)}{\sqrt{nm}}\V.
\end{align}
The main idea is to treat separately the cases $d<(\log q)^{0.05}$
and $d\ge (\log q)^{0.05}$. Note that in the main term of the
theorem, the factor $P_q(1)\psi(q)\gg (\log\log q)^{-6}$, so we need
to bound any error terms by $q\log q (\log\log q)^{-7}$.

\subsection{Small divisors}
\begin{lemma}\label{lemmaofStefanicki}
We have for some $\theta>0$,
\begin{align}
\label{offdiag} \sum_{\substack{n\neq m\\ n\equiv m \bmod \ell\\
(nm,q)=1}}
\frac{a_f(n)a_f(m)}{\sqrt{nm}}V\Big(\frac{nm}{q^2}\Big)\ll
\phi(\ell)^{-1}q (\log q)^{1-0.065}+ q^{\epsilon} \ell^{-\theta}.
\end{align}
\end{lemma}
\noindent The above lemma is a combination of Lemmas 2, 3 and 7 of
\cite{stef}. For $n$ and $m$ both close to $\ell$, a shifted
convolution problem is solved to get a power saving bound. For $n$
or $m$ much larger than $\ell$, the sum is bounded absolutely by
appealing to Shiu's \cite{shiu} estimate for sums of multiplicative
functions in arithmetic progression combined with Rankin's
\cite{rank} result that $|a_f(n)|$ is bounded by $(\log n)^{-0.065}$
on average. This gives only a $\log$ power saving. The remaining
ranges of $n$ and $m$ are trivially bounded with power savings. Thus
we have
\begin{align}
\label{arith} &\sum_{\substack{d|q\\d<(\log q)^{0.05}}} \mu(d)
\phi(q/d) \sum_{\substack{n\equiv m
\bmod q/d\\n\neq m\\(mn,q)=1}} \frac{a_f(n)a_f(m)}{\sqrt{nm}}\V\\
&\nonumber \ll \sum_{\substack{d|q\\d<(\log q)^{0.05}}} |\mu(d)|
\phi\big(\tfrac{q}{d}\big)\Big( \phi\big(\tfrac{q}{d}\big)^{-1} q
(\log q)^{0.935}+ q^{\epsilon} \big(\tfrac{q}{d}\big)^{-\theta}
\Big)\ll q (\log q)^{0.985}.
\end{align}
Now we turn to the `diagonal' terms; those with $n=m$:
\begin{align}
\label{diag} \sum_{\substack{d|q\\d<(\log q)^{0.05}}}
\mu(d)\phi\Big(\frac{q}{d}\Big) \sum_{\substack{n\\(n,q)=1}}
\frac{a_f(n)^2}{n}V\Big(\frac{n^2}{q^2}\Big).
\end{align}
First note that by Lemma \ref{condition} we have
\begin{align}
\label{diag0} \sum_{\substack{d|q\\d<(\log q)^{0.05}}}
\mu(d)\phi\Big(\frac{q}{d}\Big)&=\sum_{d|q}
\mu(d)\phi\Big(\frac{q}{d}\Big)+O(q(\log\log q)^{-9})\\
\nonumber  &=q\psi(q)+O(q(\log\log q)^{-9}).
\end{align}
The last line follows by putting $n=m=1$ in Lemma \ref{ortho}. The
sum over $n$ can be written as an integral involving the
Rankin-Selberg $L$-function $L(f\otimes f, s)$ to get
\begin{lemma}
\begin{align}
\label{diag2}\sum_{\substack{n\\(n,q)=1}}
\frac{a_f(n)^2}{n}V\Big(\frac{n^2}{q^2}\Big)= K P_q(1) \log q + K_1
P_q(1) +K_2 P_q'(1)+O(q^{-1/5}),
\end{align}
for some constants $K_i$ and $K=\text{residue}_{s=1} 2L(f\otimes
f,s)=\tfrac{6(4\pi)^k}{\pi\Gamma(k)}\|f\|^2$.
\end{lemma}
\noindent This is shown in section 4 of \cite{stef}, where is it is
also observed that $P_q'(1)\ll \nu(q) P_q(1)$. Now as $\nu(q)\ll
\frac{\log q}{\log \log q}$, we find that the diagonal terms give
the main term of Theorem \ref{mainthm}.

\subsection{Large divisors}
For $d|q$ we can express the condition $n\equiv m \bmod q/d$ using
the Dirichlet characters modulo $q/d$. We have
\begin{align}
&\phi(q/d)\sum_{\substack{n\equiv m \bmod q/d\\(mn,q)=1}}
\frac{a_f(n)a_f(m)}{\sqrt{nm}}\V\\
&\nonumber =\sum_{\substack{n,m\\(mn,q)=1}}\sum_{\chi \bmod q/d}
\frac{a_f(n)a_f(m)\chi(n)\overline{\chi}(m)}{\sqrt{nm}}\V.
\end{align}
Now the last line equals
\begin{multline}\label{integ}
\frac{1}{\pi i}\int_{(c)}\sum_{\chi \bmod q/d}
L(\f,s+\half)L(f\times
\overline{\chi},s+\half)\\
\prod_{p|d}\big(1-\tfrac{a_f(p)\chi(p)}{p^{s+1/2}}+\tfrac{\chi(p^2)}{p^{2s+1}}\big)\prod_{p|d}\big(1-\tfrac{a_f(p)\overline{\chi}(p)}{p^{s+1/2}}+\tfrac{\overline{\chi}(p^2)}{p^{2s+1}}\big)
G(s)q^{2s}\frac{ds}{s}.
\end{multline}
Here the extra Euler products arise due to the condition $(nm,q)=1$.
On taking $c=1/\log q$ and using Lemma \ref{lemmabound} from the
final section, we have that (\ref{integ}) is less than a constant
multiple of
\begin{align}
\frac{q}{d} \log q \log \log q
\prod_{p|d}\big(1+\tfrac{10}{\sqrt{p}}\big) \int_{\frac{1}{\log
q}}\Big| (s+1)^4 G(s)q^{2s}\frac{ds}{s}\Big|.
\end{align}
The last integral is easily found to be $\ll \log \log q$. Thus we have
\begin{multline}
\sum_{\substack{d|q\\ d\ge (\log q)^{0.05}}} \mu(d) \phi(q/d)
\sum_{\substack{n\equiv m \bmod q/d\\(mn,q)=1}}
\frac{a_f(n)a_f(m)}{\sqrt{nm}}\V\\
\ll q\log q (\log \log q)^2 \sum_{\substack{d|q\\ d\ge (\log
q)^{0.05}}} \frac{|\mu(d)| \prod_{p|d}(1+\tfrac{10}{\sqrt{p}}) }{d}.
\end{multline}
This is $\ll q\log q (\log \log q)^{-7}$ by Lemma \ref{condition}.

\section{An upper bound for the second moment}

In this section we will need a different form of the approximate
functional equation.
\begin{lemma}[\bf Approximate functional equation]
For a primitive character $\chi$ modulo $q$ and $\Re(s)=0$ we have
\begin{multline} \label{approx2}L(\f,\half+s)=\sum_n
\frac{a_f(n)\chi(n)}{n^{\frac{1}{2}+s}}W_s\Big(\frac{n}{q}\Big)\\
+\imath_{\chi}\frac{(2\pi)^{2s}\Gamma(k/2-s)}{q^{2s}\Gamma(k/2+s)}
\sum_n
\frac{a_f(n)\overline{\chi}(n)}{n^{\frac{1}{2}-s}}W_{-s}\Big(\frac{n}{q}\Big),
\end{multline}
where $W_s(x)=\frac{1}{2\pi i}\int_{(c)} \frac{\Gamma(k/2+s+y)}{
\Gamma(k/2+s)}e^{y^2}(2\pi x)^{-y}\frac{dy}{y}$ for any $c>0$. We
have that $W_s(x) \ll_C \min\{ 1, |s+1|^C x^{-C} \}$ for any $C>0$,
so that the sums above are essentially supported on $n<
q^{1+\epsilon}$ for any $\epsilon>0$.
\end{lemma}
\noindent This is a special case of Theorem 5.3 in \cite{iwakow}.
The bound on $W_s(x)$ can be seen by moving the line of integration
far to the right when $x>1$ and just to the left of the origin when
$x\le 1$.

We show the following.

\begin{lemma}\label{lemmabound}
For $\Re(s)\ge 0$ and integers $q$ and $r$ we have
\begin{align}
\label{lemmab} \Big| &\sum_{\chiq} L(f\times \chi,\half+s)L(f\times
\overline{\chi},\half+s)\\
&\nonumber \hspace{1in} \prod_{p|r}\big(1-\tfrac{a_f(p)\chi(p)}{p^{s+1/2}}+\tfrac{\chi(p^2)}{p^{2s+1}}\big)\big(1-\tfrac{a_f(p)\overline{\chi}(p)}{p^{s+1/2}}+\tfrac{\overline{\chi}(p^2)}{p^{2s+1}}\big) \Big| \\
&\nonumber \ll |s+1|^{4} q \log q \log \log q \prod_{p|r}
\big(1+\tfrac{10}{\sqrt{p}}\big).
\end{align}
\end{lemma}
\proof The statement of the lemma certainly holds for $\Re(s)\ge 2$,
so by the Phragmen-Lindel\"{o}f principle (cf. \cite[Theorem
5.53]{iwakow}) it is enough to prove the result for $\Re(s)=0$. Thus
throughout this proof, assume that $s$ is purely imaginary. First we
demonstrate the bound
\begin{align}
\label{primbound} \summ_{\chid} |L(f\times \chi,\half+s)|^2\ll
|s+1|^{4} d\log d.
\end{align}
By Lemma \ref{approx2} it is enough to bound
\begin{align}
\label{cs} &\summ_{\chid} \Big|\sum_n
\frac{a_f(n)\chi(n)}{n^{\frac{1}{2}+s}}W_s\Big(\frac{n}{d}\Big)\Big|^2\\
&\nonumber\ll \sum_{\chid} \Big|\sum_{n\le d}
\frac{a_f(n)\chi(n)}{n^{\frac{1}{2}+s}}W_s\Big(\frac{n}{d}\Big)\Big|^2+\sum_{\chid}
\Big|\sum_{n>d}
\frac{a_f(n)\chi(n)}{n^{\frac{1}{2}+s}}W_s\Big(\frac{n}{d}\Big)\Big|^2.
\end{align}
By a large sieve inequality (cf. \cite[section 7.5]{iwakow}) and the
Rankin-Selberg bound $\sum_{n \le x} a_f(n)^2 \ll x$, we have
\begin{align}
\sum_{\chid} \Big|\sum_{n\le d}
\frac{a_f(n)\chi(n)}{n^{\frac{1}{2}+s}}W_s\Big(\frac{n}{d}\Big)\Big|^2\ll
 d\sum_{n\le d} \frac{a_f(n)^2}{n}\ll  d\log d.
\end{align}
The remaining sum in (\ref{cs}) equals, by orthogonality,
\begin{align}
\phi(d)\sum_{\substack{1\le r \le d\\ (r,d)=1}} \Big|
\sum_{\substack{n>d \\n\equiv r \bmod d}}
\frac{a_f(n)}{n^{\frac{1}{2}+s}}W_s\Big(\frac{n}{d}\Big)\Big|^2.
\end{align}
Dropping the condition $(r,d)=1$, this is bounded by
\begin{align}
&\label{lem1} d\sum_{h\ge 0} \sum_{ n>d}
\frac{a_f(n)a_f(n+hd)}{n^{\frac{1}{2}+s}
(n+hd)^{\frac{1}{2}-s}}W_s\Big(\frac{n}{d}\Big)\overline{W_s\Big(\frac{n+hd}{d}\Big)}\\
&\ll \nonumber d  \sum_{h\ge 0} \sum_{ n>d}
\frac{a_f(n)^2+a_f(n+hd)^2}{n^{\frac{1}{2}}
d^{\frac{1}{2}}}\Big(\frac{d}{n}\Big)(1+h)^{-3}|1+s|^4 \ll d
|1+s|^4.
\end{align}
This establishes (\ref{primbound}).

Now, the left hand side of (\ref{lemmab}) bounded by
\begin{align}
\prod_{p|r} \big(1+\tfrac{10}{\sqrt{p}}\big)\sum_{\chiq} |L(f\times
\chi,\half+s)|^2.
\end{align}
We have
\begin{align}
\label{lemmaproof}   \sum_{\chiq} |L(f\times &\chi,\half+s)|^2=  |L(f\times \chi_0,\half+s)|^2 \\
&\nonumber+ \sum_{d|q} \summ_{\chid} |L(f\times \chi,\half+s)|^2
\prod_{p|\frac{q}{d}}\Big|1-\frac{a_f(p)\chi(p)}{p^{\frac{1}{2}+s}}+\frac{\chi(p^2)}{p^{2s}}\Big|^2.
\end{align}
Using (\ref{primbound}) we have
\begin{align}
\sum_{\chiq} |L(f\times \chi,\half+s)|^2 &\ll |s+1|^{4} \sum_{d|q} d
\log d
\prod_{p|\frac{q}{d}}(1+\tfrac{10}{\sqrt{p}})\\
&\nonumber\ll |s+1|^{4}q\log q\prod_{p|q}(1+\tfrac{1}{p})\\
&\nonumber\ll |s+1|^{4}q\log q \log \log q.
\end{align}
\endproof

\vspace{.3in}

{\bf Acknowledgements.} We would like to thank Prof. K.
Soundararajan for his encouragement. We also thank the referee for a
careful reading and helpful suggestions on exposition. The first
author is supported by a research fellowship from an Academic
Research Fund Tier 1 grant at Nanyang Technological University. The
second author is partially supported by the National Science
Foundation. The third author is financed by the ANR project
``Aspects Arithm\'etiques des Matrices Al\'{e}atoires et du Chaos
Quantique''.

\bibliographystyle{amsplain}

\bibliography{AA5891_references}

\end{document}